\documentclass[10pt]{article}
\newcommand{\documentdate}{24 November 2017}

\usepackage{latexsym,graphicx,color}
\usepackage[a4paper]{geometry}

\newcommand{\calB}{{\cal B}}
\newcommand{\calP}{{\cal P}}
\newcommand{\calQ}{{\cal Q}} 
\newcommand{\calS}{{\cal S}}
\newcommand{\hatq}{\hat{q}}
\newcommand{\hatt}{\hat{t}}
\newcommand{\tim}[1]{\;\; \mbox{#1} \;\;}

\newcommand{\beqn}[1]{\begin{equation}\label{#1}}
\newcommand{\eeqn}{\end{equation}}
\newcommand{\req}[1]{(\ref{#1})}
\newcommand{\eqdef}{\stackrel{\rm def}{=}}
\newcommand{\bigint}{\displaystyle \int}

\title{A note on using performance and data profiles \\ for training algorithms}

\author{M. Porcelli\thanks{%
Universit\`a degli Studi di Firenze, Dipartimento di Ingegneria Industriale, 
viale G.B. Morgagni 40, 50134, Firenze, Italy
Email: margherita.porcelli@unifi.it}\,
~and  Ph. L. Toint\thanks{%
Namur Center for Complex Systems (NAXYS),
University of Namur,
61, rue de Bruxelles, B-5000 Namur, Belgium.
Email: philippe.toint@unamur.be}}
\date{\documentdate}

\begin{document}

\maketitle

\begin{abstract}
It is shown how to use the performance and data profile benchmarking tools to
improve algorithms' performance.  An illustration for the BFO derivative-free
optimizer suggests that the obtained gains are potentially significant.
\end{abstract}
    {\bf Keywords:} algorithmic design, algorithms' training, trainable codes,
    derivative-free optimization.\\

\noindent
{\bf Mathematics Subject Classification:} 65K05, 90C56, 90C90.

\section{Introduction}

Making algorithms efficient and reliable is obviously desirable for both their
designers and their users.  Since most algorithms involve parameters, it is
therefore important to choose them well. In an attempt to do so, the authors
\cite{PorcToin17} proposed BFO, a derivative-free optimization algorithm which
is trainable, in the sense that it contains an internal procedure to select
its algorithmic parameters to improve algorithmic performance, both from the
point of view of the designer (using a large collection of diverse
benchmarking cases) and of the user (focusing on a possibly more specific
class of applications). Moreover BFO is also designed so that it can be used to
train other codes. Obviously, improving performance requires a workable
definition of this concept.  As in \cite{PorcToin17}, we assume that
performance of an algorithm on a given problem can be measured by a number and
that better performance corresponds to smaller such numbers. To make things
concrete, and since we will be concerned below with derivative-free
optimization, we shall consider from now on that performance is given by the
number of objective function evaluations required by a solver to solve a given
optimization problem\footnote{In other contexts, possible measures involve
error on the result, computing time, memory usage, etc.} Given a vector of
algorithmic parameters $q$ and a collection of benchmarking problems $\calP$,
algorithmic performance was then measured using one of two classical
techniques. The first is the 'average' performance over all test problems, the
second, inspired by robust optimization, is the average of the worst
performance obtained for a slight variation of $q$. In the first case, one
then attempts to improve algorithmic performance by approximately minimizing
the \emph{average training function}
\cite{AudeOrba06,AudeDangOrba10,AudeDangOrba14}
\begin{equation}\label{AO}
\min_{q}\ \phi^A_\calP(q)   
\end{equation}
where $\phi^A_\calP(q)$ counts the total number of evaluations of the problem
objective function to solve all the problems in  $\calP$. In the second case,
improvement is sought by approximately minimizing the \emph{robust training function}
\begin{equation}\label{RO}
  \min_{q} \  \phi^R_\calP( q )
  \tim{ where } \phi^R_\calP( q ) \eqdef \max_{\hatq \in \calB(q) } \  \phi^A_\calP(\hatq) 
\end{equation}
with $\calB(q)$ being a local box centered at $q$ allowing perturbations of 
each algorithmic parameter. In both cases, the minimization of the
training function also typically involves bound constraints on the admissible
range of each algorithmic parameter. It was shown in \cite{PorcToin17} that an
approximate local minimization of either of these training functions can bring
substantial improvements in efficiency and reliability. The final comparison
(and that with other derivative-free approaches) was then reported using 
the now widely accepted performance and data profile techniques (see
\cite{DolaMoreMuns06} for the first and \cite{MoreWild09} for the second).

The purpose of the present short note is to explain how it is possible to use these
two latter benchmarking measures directly for training, instead of merely for comparison.
As in \cite{PorcToin17}, we focus on the BFO derivative-free solver because it
directly implements the relevant tools, but we stress that the approach is not
limited to this particular case.

The paper is organized as follows. We first briefly recall, in Section~\ref{prof-s},
the definition of performance and data profiles given in
\cite{DolaMoreMuns06,MoreWild09} and then derive the new training measures and
associated training procedures in Section~\ref{trainmeas-s}. A numerical illustration
is reported in Section~\ref{experiments-s}.

\section{Performance and data profiles}\label{prof-s}
  
Let $\calS$ be a set of solvers (or solver variants) and let $\calP$ be a set
of benchmarking problems of cardinality $|\calP|$.  {\em Performance profiles}
are defined in terms of a performance measure $t_{p,s} > 0$ obtained for each
$p\in \calP$ and $s\in \calS$. As above, we will consider here that $t_{p,s} >
0$ is the number of function evaluations required to satisfy a user-defined
convergence test.  For each $p\in \calP$, let $\hatt_{p,s} = \min_{s\in \calS}
t_{p,s}$ and define $r_{p,s} = t_{p,s}/ \hat t_{p,s} $ to be the performance
ratio, so that the best solver $s$ for a particular problem $p$ attains the
lower bound $r_{p,s} = 1$. We set $r_{p,s} = \infty$ when solver $s$ fails to
satisfy the convergence test on problem $p$. For $\tau \ge 1$, each solver
$s\in \calS$ and each problem $p\in \calP$, one then defines
\[
 k(r_{p,s}, \tau) = \left \{ \begin{array}{ll}
                              1 & \mbox{ if } r_{p,s} \le \tau, \\
                              0 & \mbox{ otherwise. }
                             \end{array} \right. 
\]
The performance profile for solver $s$ is then given by the function
\[
 p_s(\tau) = \frac{1}{|\calP|}\sum_{p \in \calP}k(r_{p,s}, \tau), 
 \quad \tau \ge 1.
\]
By definition of $t_{p,s}$, $p_s(1)$ is the fraction of problems 
for which solver $s$ performs the best, $p_s(2)$ gives the fraction of problems
for which the solver's performance is within a factor of 2 of the best, and that 
for $\tau$ sufficiently large, $p_s(\tau)$ is the fraction of problems solved by $s$.
More generally, $p_s(\tau)$ can be interpreted as the probability for solver $s\in S$
that the performance ratio $r_{p,s}$ is within a factor $\tau$ of the best possible ratio.
Therefore, $p_s(1)$ measures efficiency of the solver while its robustness
(high probability of success on the set $\calP$) is measured in terms of 
$p_s(\infty)$. A key feature of performance profiles is that they give information
on the \emph{relative} performance of several solvers \cite{DolaMore02,MoreWild09},
which therefore strongly depends of the considered set $\calS$ of competing
solvers or algorithmic variants \cite{GoulScot16}.

In order to provide a benchmarking tool that gives the behaviour of a solver 
independently of the other solvers in $\calS$, Mor\'e and Wild
\cite{MoreWild09} proposed the {\em data profile} measure 
motivated by the user interest in the percentage of problems that can be
solved with  a certain computational ``budget''.
For $\nu > 0$ and each $s\in \calS, p\in \calP$, one defines
\[
 g(t_{p,s}, \tau) = \left \{ \begin{array}{ll}
                              1 & \mbox{ if } t_{p,s} \le \nu (n_p + 1), \\
                              0 & \mbox{ otherwise, }
                             \end{array} \right.                             
\]
where $n_p$ is  the number of variables in $p\in \calP$. The scaling by
$n_p+1$ is intended to consider the computational budget as in ``simplex
gradient'' evaluations, rather than directly in function evaluations.
The data profile for solver $s\in \calS$ is then given by
\[
 d_s(\nu) = \frac{1}{|\calP|}\sum_{p \in \calP}g(t_{p,s}, \nu), 
 \quad \nu >0,
\]
and measures the percentage of problems that can be solved with $\nu$
``simplex gradient'' evaluations.

\section{New training measures and how to use them}\label{trainmeas-s}

We observe that, by definition, the plots of the performance and data 
profiles are staircase graphs and that, by the above discussion, the higher
the curve corresponding to a solver, the better is its performance.
This trivial observation suggests two new training strategies that simply
consists in finding the parameter configuration that maximize the area
under the staircase graph generated by the performance or data profiles,
respectively.

Let $\calQ$ be the set of acceptable algorithmic parameters, $q \in \calQ$ be
a parameter configuration and let $s_q$ be the solver variant with parameter
configuration $q$.  Consider data profiles first. We can define for each $q\in
\calQ$ the \emph{data profile training function}
\[
\phi^D_\calP(q) \eqdef \int_{\nu_{\min}}^{\nu_{max}} d_{s_q}(\nu) \,d\nu,
\]
where $0 \le \nu_{\min} < \nu_{max}$ are user-specified values identifying a
'range of computational budgets' of interest, and then consider the
corresponding data profile training problem
\begin{equation}\label{DP}
\max_{s \in \calQ}\  \phi^D_\calP(q).
\end{equation}
%We observe that problem \req{DP} can be interpreted as the maximization 
%of a data profile ``density''.

The analogous problem for performance profiles is less obvious since, as
discussed above, the computation of $p_{s_q}(\nu)$ depends on the behaviour of
more than one solver, that is, in our case, on the performance of the trained
solver with respect to different values of its algorithmic parameters $q$.  We
therefore propose to proceed sequentially from an initial parameter configuration
indexed by $i=0$ and to evaluate the performance for a particular $q$ by
always comparing it to that obtained for $q_0$.  Given the profile window
$[\tau_{\min}, \tau_{max}]$ for some $1\le \tau_{\min} < \tau_{max} $ and the
initial algorithmic configuration $q_0\in \calQ$, we define the
\emph{performance profile training function}
$\phi^P_\calP$ by
\beqn{PPobjf}
%\begin{array}{lcl}
\phi^P_\calP( q ) \eqdef \bigint_{\tau_{\min}}^{\tau_{max}} 
 \left[p_{s_q}(\tau)-p_{s_{q_0}}(\tau)\right] \,d\tau.
%\end{array}
\eeqn
Training then corresponds to solving (possibly very approximately) the
performance profile training problem
\begin{equation}\label{PP}
\max_{q \in \calQ}\  \phi^P_\calP( q ).
\end{equation}

In order to evaluate $\phi^P_\calP$ and $\phi^D_\calP$ in \req{DP}-\req{PP}
respectively, one has to provide enough information to compute the profiles
$p_s(\tau)$  and $d_s(\nu)$ during the training optimization process.

Let $q \in \mathcal{Q}$ be a parameter configuration and let $s_q$
be the (BFO) algorithmic variant using parameters $q$. Let the profiles
windows $[\tau_{min}, \tau_{max}]$ and $[\nu_{\min}, \nu_{max}]$ be given.
We compare different parameter configurations declaring that the problem
$p$ with objective function $f_p$ is solved by the variant $s_q$ as soon as it
produces an approximate solution $x_q$ such that
\begin{equation}\label{convtest}
f_p(x_q) \le f_p^* + \chi( f_p(\bar x)- f_p^*) \eqdef c_p
\end{equation}
where $\bar x$ is the starting point for the problem $p$, $f_p^*$ is an
approximation of the smallest obtainable value of $f_p$ and $\chi\in[0,1]$ is
a tolerance. The test \req{convtest} therefore compares the function value
reduction $f(\bar x)-f(x_q)$ achieved by $x_q$ relative to the best possible
reduction $f(\bar x) - f_*$ \cite{MoreWild09}. We say that $c_p$, as defined
in \req{convtest}, is the \emph{cut-off value} for problem $p$.

Given an initial parameter configuration $q_0$, a starting point $\bar x$ and
a tolerance $\chi >0$, the training strategy proceed as follows.  First,
starting from $\bar x$, the solver variant $s_{q_0}$ is run over the set
$\mathcal{P}$ with high accuracy in order to evaluate the best objective found
$f_p^*$ for each $p \in \mathcal{P}$ and the resulting cut-off value $c_p$.
Then, the number of function evaluations needed to the solver variant
$s_{q_0}$ to reach $c_p$, that is the value $t_{p, s_{q_0}}$, is retrieved.
If data-profile training is considered, this enough to compute the
corresponding value of the objective $\phi^D_\calP(q_0)$. The initial
objective function value for performance-profile training is initialized to
zero (see \req{PPobjf}).  Optimizing the relevant objective function
(i.e.\ \req{DP} or \req{PP}) can then be conducted (using BFO with its default
parameters and its standard termination test in our case), in the course of
which the solver variant is run again with better and better values of the
algorithmic parameters $q$, the performance measures $t_{p,s_q}$
being always computed with respect to the initial cut-off value $c_p$.

\section{Numerical illustration}\label{experiments-s}

We now illustrate the above proposals by reporting some results obtained when
training the BFO derivative-free optimization package by modifying its
internal algorithmic parameters.

\subsection{Experimental setup}

The inner details of the BFO method are of little interest here (we refer the
interested reader to \cite{PorcToin17} for a full description).  It is
enough for our present purposes to describe it as a direct-search optimizer
evaluating the relevant objective function at points on a randomly oriented
variable meshsize grid (in a process called the poll-step) and accepting an
improved function value whenever it satisfies a 'sufficient decrease'
condition relative to the current grid meshsize. So-called 'inertia direction'
are also computed using a number of past iterates and are priviledged when
constructing the grid. The minimization is terminated when the grid meshsize
becomes smaller than a user-supplied threshold $\epsilon$. The BFO algorithmic
parameters considered for training in our present experiments are presented in
Table \ref{partab}.

\begin{table}[ht]
{\small 
\centering
\begin{tabular}{l|c|l}
  \hline
        Parameters & Type  & Description \\
  \hline
  $\alpha$ & $c$ &The grid expansion factor\\
  $\beta$ & $c$ &The grid shrinking factor\\
  $\gamma$ & $c$ & The maximum grid expansion factor\\
  $\delta$ & $c$ &The initial stepsize vector \\
  $\eta$ &  $c$ &The sufficient decrease fraction in the poll step \\
  $\tt inertia$ &  $i$ &The number of iterations for the inertia direction  \\
 \hline
\end{tabular}  
\caption{BFO parameters selected for training.}\label{partab}} \end{table}

We define $\calP$ to be the set of benchmarking problems used in
\cite{PorcToin17} and consisting in 55 bound-constrained problems with
continuous variables of small dimensions extracted from the {\sf CUTEst}
library \cite{GoulOrbaToin15b}. The list of problem names with their dimension
is given in Table \ref{testpb}. The solution of each test problem is attempted
setting $\epsilon = 10^{-12}$ in the  BFO convergence test and allowing 10000
function evaluations at most. 

\begin{table}[ht]
\centering
{\footnotesize
\begin{tabular}{l c l c lc lc lc }
  \hline
   Name &$n$&   Name &$n$&   Name &$n$&   Name &$n$&   Name &$n$ \\
ALLINIT &4 &	HADAMALS  &4&	HS38	&4&	MDHOLE	&2&	PENTDI	&5   \\
BDEXP	&10&	HARKERP2  &10&	HS3     &2&	NCVXBQP1&10&	POWELLBC&6\\
BIGGSB1	&10&	HART6     &6&	HS3MOD  &2&	NCVXBQP2&10&	PROBPENL&10\\
CAMEL6	&2&	HATFLDA   &4&	HS45    &5&	NCVXBQP3&10&	PSPDOC	&4\\
CHARDIS0&10&	HATFLDB   &4&	HS4     &2&	NONSCOMP&10&	QUDLIN	&12\\
CHEBYQAD&4&	HATFLDC   &9&	HS5     &2&	OSLBQP	&8&	S368	&8\\
CVXBQP1&10&	HIMMELP1  &2&	KOEBHELB&4&	PALMER1A&6&	SIMBQP	&2\\
EG1	&3&     HS110     &10&	LINVERSE&9&	PALMER2B&4&	SINEALI	&4\\
EXPLIN	&12&    HS1       &2&	LOGROS  &2&	PALMER3E&8&	SPECAN	&9\\
EXPLIN2	&12&	HS25      &3&	MAXLIKA &8&     PALMER4A&6&	WEEDS	&3\\
EXPQUAD	&12&	HS2       &2&	MCCORMCK&10&	PALMER4	&4&	YFIT	&3  \\
  \hline
\end{tabular} 
\caption{The benchmark problem set $\calP$: name and dimension $n$.} \label{testpb}}
\end{table}

Starting from the initial parameter configuration $q_0$ in Table \ref{parset},
the four optimization problems \req{AO}-\req{RO}-\req{DP}-\req{PP} are
approximately solved imposing bound constraints on the parameters with bounds
$l$ and $u$ reported in Table \ref{parset}.  The local box $\calB( q)$ in
\req{RO} is defined for continuous parameters as the Cartesian product of the
intervals $\left[0.95\, q, 1.05\, q \right]$ for continuous parameters and
$\{q\}$ for discrete ones, allowing perturbations of each
continuous algorithmic parameter by at most $5\%$.

\begin{table}[ht]
{\small 
\centering
\begin{tabular}{l c c c c  c c r}
  \hline
   & $\alpha$ & $\beta$ & $\gamma$ & $\delta$ & $\eta$ & $ \tt inertia$ & \\		
  \hline
      $q_0$ &   1.5   &   1/3   & 5    &  1     &   $10^{-1}$ & 10   &  \\
      $l$ &    1    &  0.01   & 1    &  0.25  &  $10^{-4}$    &  5   &  \\
      $u$ &    2    &  0.95   & 10   &  10    &   0.5         & 30   &  \\
     \hline 
  \end{tabular}  
\caption{Starting parameter configuration $q_0$ and lower/upper bounds in the 
training optimization problems.} \label{parset}}
\end{table}

As in \cite{PorcToin17}, we set the BFO termination threshold $\epsilon =
10^{-2}$ when solving the training minimization problems
\req{AO}-\req{RO}-\req{DP}-\req{PP}, and $\epsilon = 10^{-1}$ for the
approximate solution of the inner minimization problem in \req{RO}.  We also
set an upper bound of 200 parameter configuration trials.
%and reset the random seed in BFO to the value 0.
Finally, the cut-of values used for data and performance profile strategies
are obtained by solving each problem with $\epsilon = 10^{-12}$ and
default parameters given in Table~\ref{parset} and the training is run using
$\chi = 10^{-4}$ in \req{convtest}. Experiments were carried out using Matlab
R2016b on Intel Core i7 CPU 920 @ 2.67GHz x 8 12GB RAM.

Even if the training process using approximate minimizations of the relevant
objective function guarantees improvements on the initial guess $q_0$, it is
important to remember that there is absolutely no guarantee of reaching
a local solution of the training problem, not to mention a global one.

\subsection{Results}

We report in Table \ref{partunecc} the values of the trained BFO parameters
obtained using the four training strategies. Values of parameters using the
profile training strategies are obtained setting the profile windows 
$[\nu_{\min}, \nu_{\max}]=[0, 2000]$ and $[\tau_{\min}, \tau_{\max}] =[1,20]$, 
for the objectives in \req{DP} and \req{PP}, respectively.
Table \ref{partunecc} also reports the gain (in percentage and as measured
with the relevant objective function) in the number of problem function
evaluations achieved by the training process.

\begin{table}[ht]
{\small 
\centering
\begin{tabular}{l c c c c  c c  r}
  \hline
   & $\alpha$ & $\beta$ & $\gamma$ & $\delta$ & $\eta$ & $ \tt inertia$  & gain\\
     \hline
       $q_{A}$         &   1.2     &   1/3   & 9.7    &  0.25   &  $10^{-4}$ & 10   & 17\%  \\
       $q_{R}$         &   2       &   0.25  & 8.5    &  0.25   &  $10^{-4}$ & 10   & 7\%  \\
       $q_{P}$  &   1.5     &   1/3   & 5.9    &  1      &  $10^{-4}$ & 11   & 1\% \\  
 % same results for PP in [1,5]. IN [15,20] It makes no sense training more.
       $q_{D}$ &   1.6    &    0.32  & 4.8    &  0.25   &  $10^{-4}$ & 11   & 7\%  \\
       
  \hline
  \end{tabular}  
\caption{Values of the trained parameters using different training
  strategies.}
\label{partunecc}}
\end{table}

These results show that potentially large gains in average number of function
evaluations may be obtained by training with the average strategy (and to a lesser
extent with the robust strategy), which is coherent with the findings of
\cite{PorcToin17}. Of course, this says little about the distribution of these
improvements across test problems, as is suggested by the fact that the
improvements are more modest in terms of performance and data profiles.

\begin{figure}[htb]
\centering
  \includegraphics[width=0.47\textwidth]{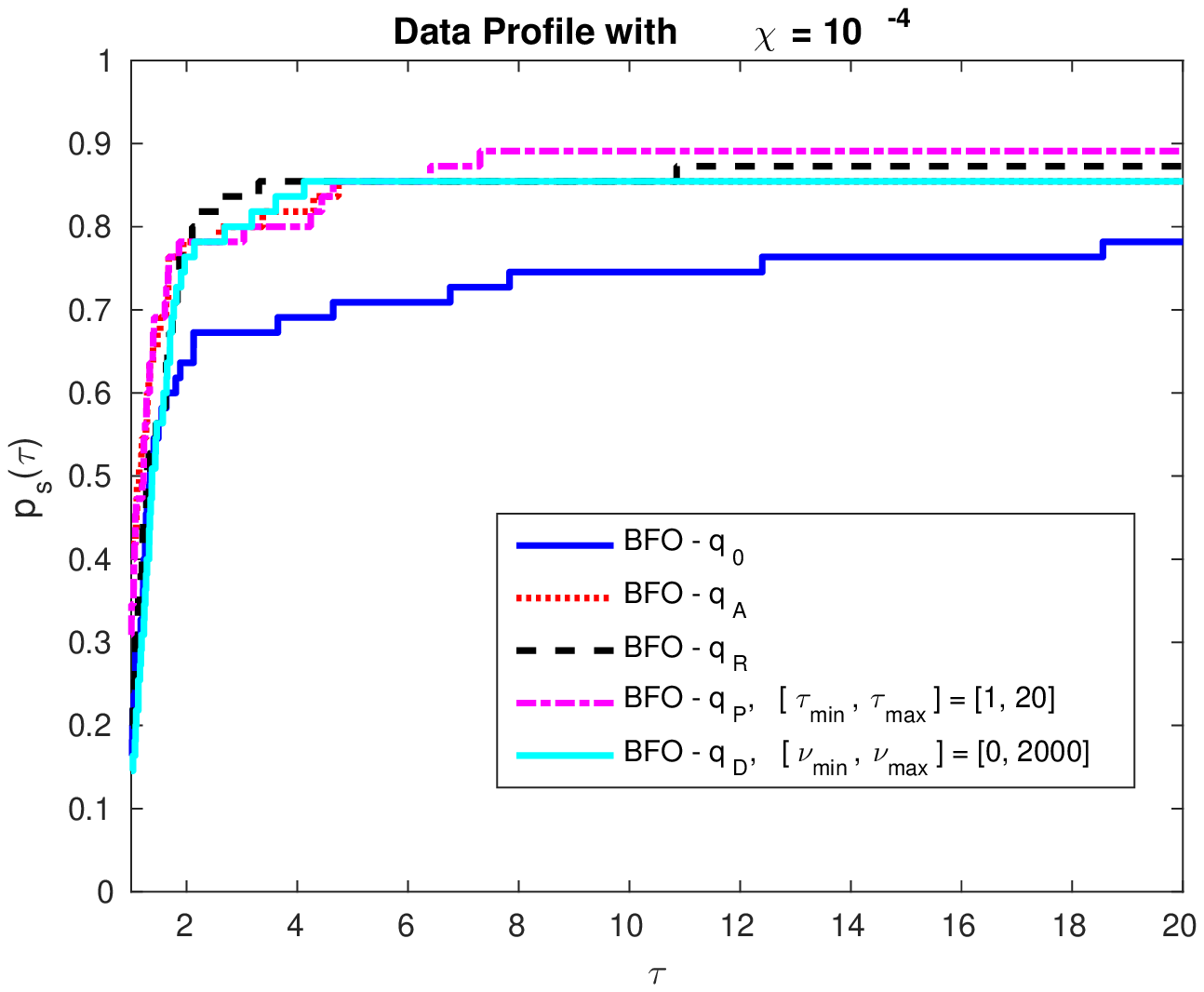}
  \includegraphics[width=0.47\textwidth]{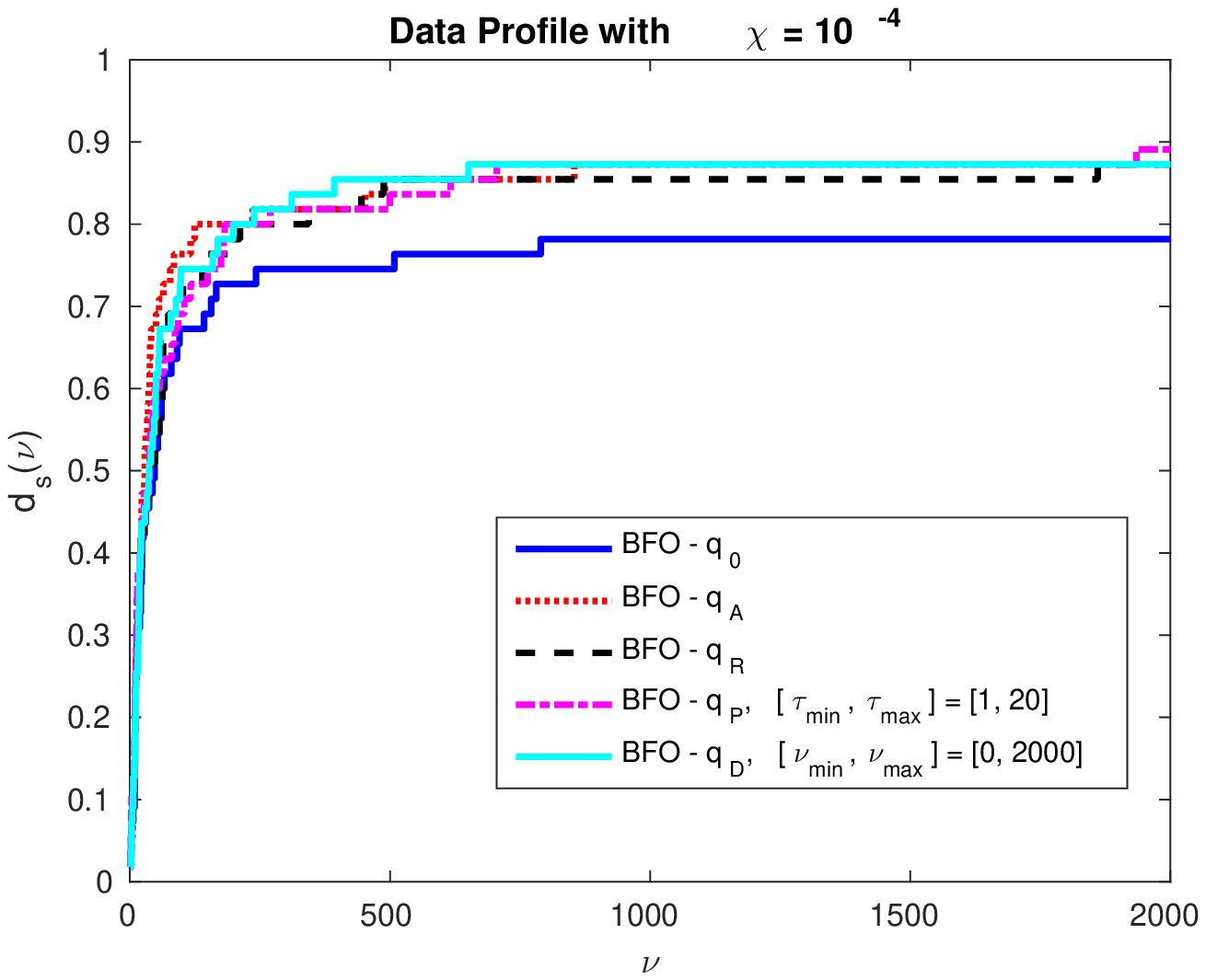}
  \caption{Performance (left) and data (right)  profiles for BFO with different 
  algorithmic parameters. Parameters $q_{P}$ and $q_{D}$ are trained in the 
  intervals $[1, 20]$ and $[0, 2000]$, respectively.} 
  \label{figure:fig1}
\end{figure}

\begin{figure}[ht]
\centering
  \includegraphics[width=0.47\textwidth]{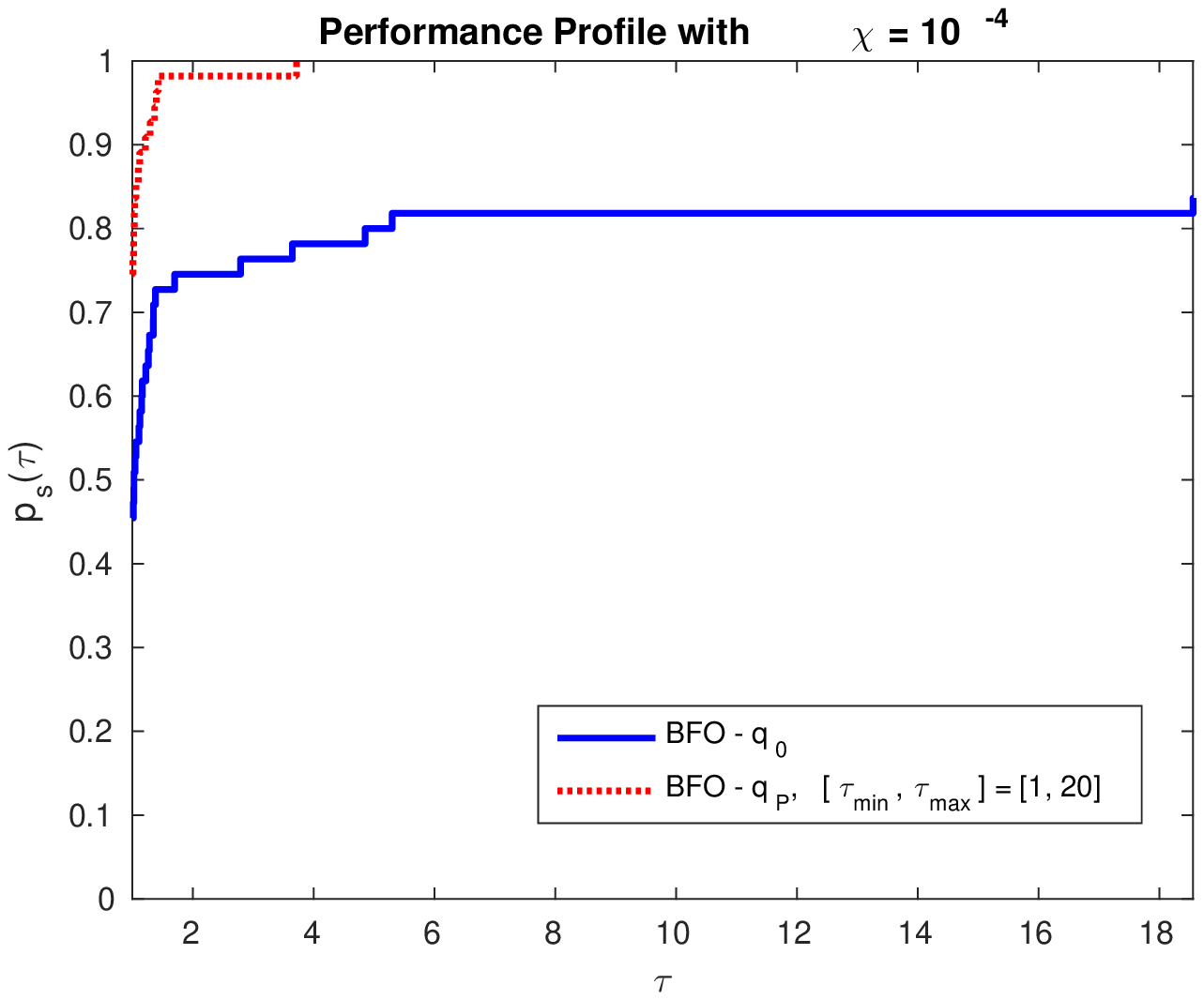}
  \includegraphics[width=0.47\textwidth]{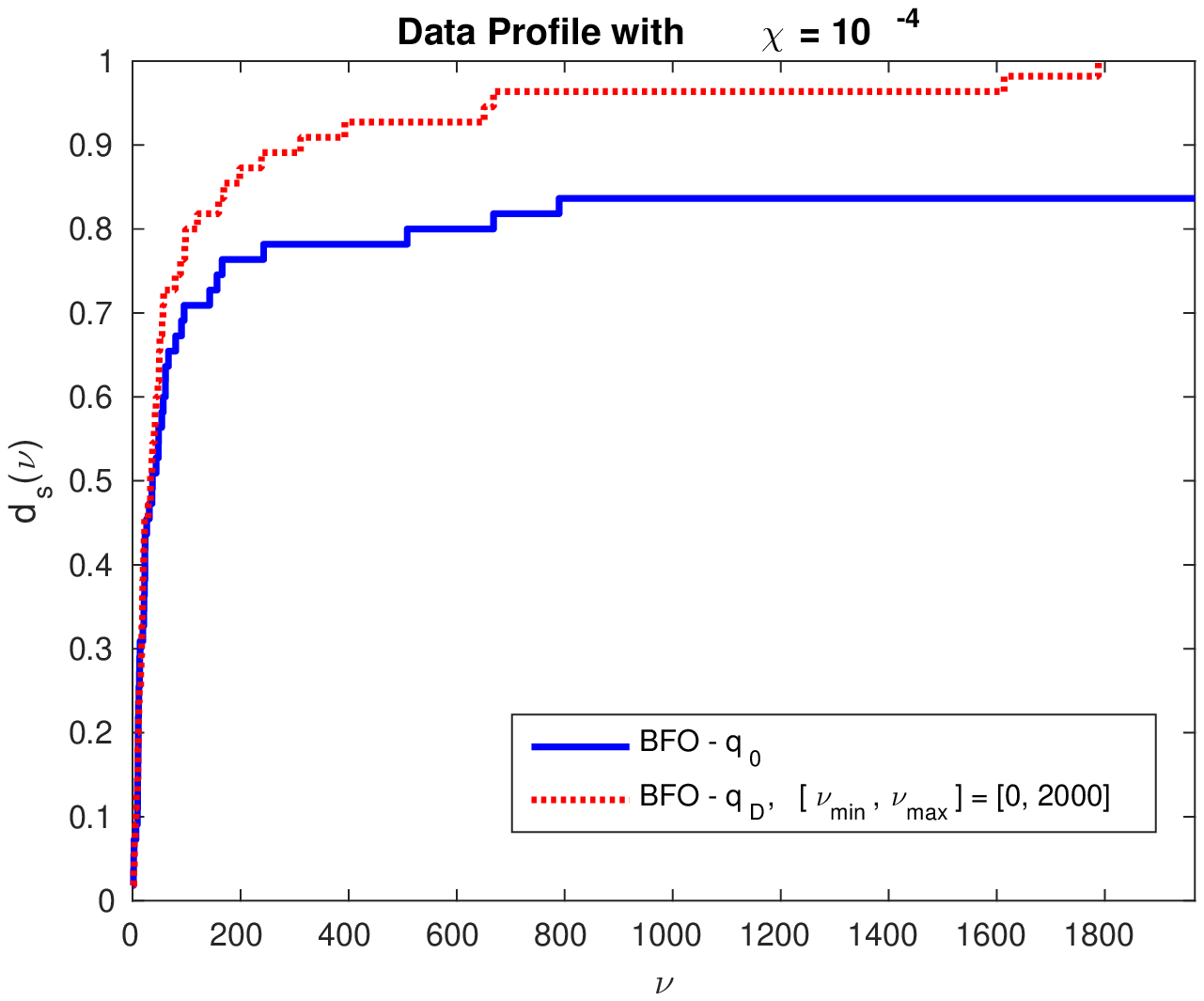}
 \caption{Performance profile for BFO with default parameters $q_0$ against BFO
 with $q_{P}$ (left), and data profiles of BFO  with default parameters $q_0$ 
 against BFO with parameters $q_{D}$ (right). Parameters $q_{P}$ and
 $q_{D}$ are trained in the default intervals $[1, 20]$ and $[0, 2000]$, respectively.} 
  \label{figure:fig2}
\end{figure}

Figure~\ref{figure:fig1} shows that, despite the lack of guarantee of global
optimality, BFO with $q_{P}$ gives the best performance in 
terms of performance profile, while  BFO with $q_{D}$ is best in terms of 
data profiles. Figure~\ref{figure:fig2} indicates that the performance
improvements are also clear when using performance or data profiles, and thus,
unsurprinsingly, that the measure of improvement reported in
Table~\ref{partunecc} for these strategies might be misleading.

We now focus on the behaviour of BFO trained by the new profile strategies
and discuss the effect on performance of varying the training windows.  From
the definitions, we would expect a profile window with small values
(i.e.\ $\tau_{\max}$ relatively modest) to boost performance, while a window
with larger values (substantial $\tau_{\min}$) to result in better
reliability.  Because the performance profile result shows little room for
improvement either in efficiency or reliability (as shown by
Figure~\ref{figure:fig2}), we illustrate these effects (and their limits) with
using data-profile training.

\begin{figure}[htb]
\centering
  \includegraphics[width=0.47\textwidth]{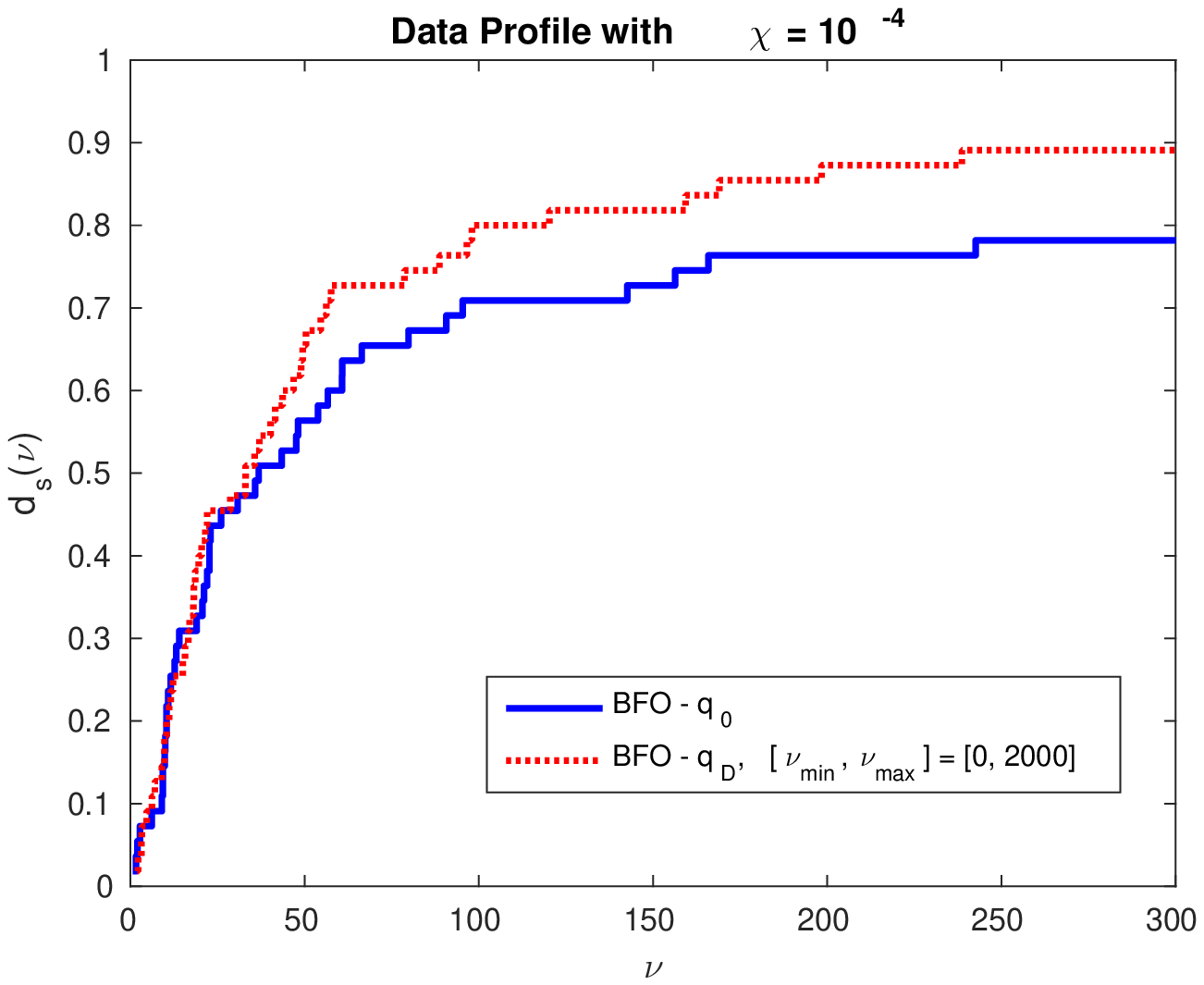}
  \includegraphics[width=0.47\textwidth]{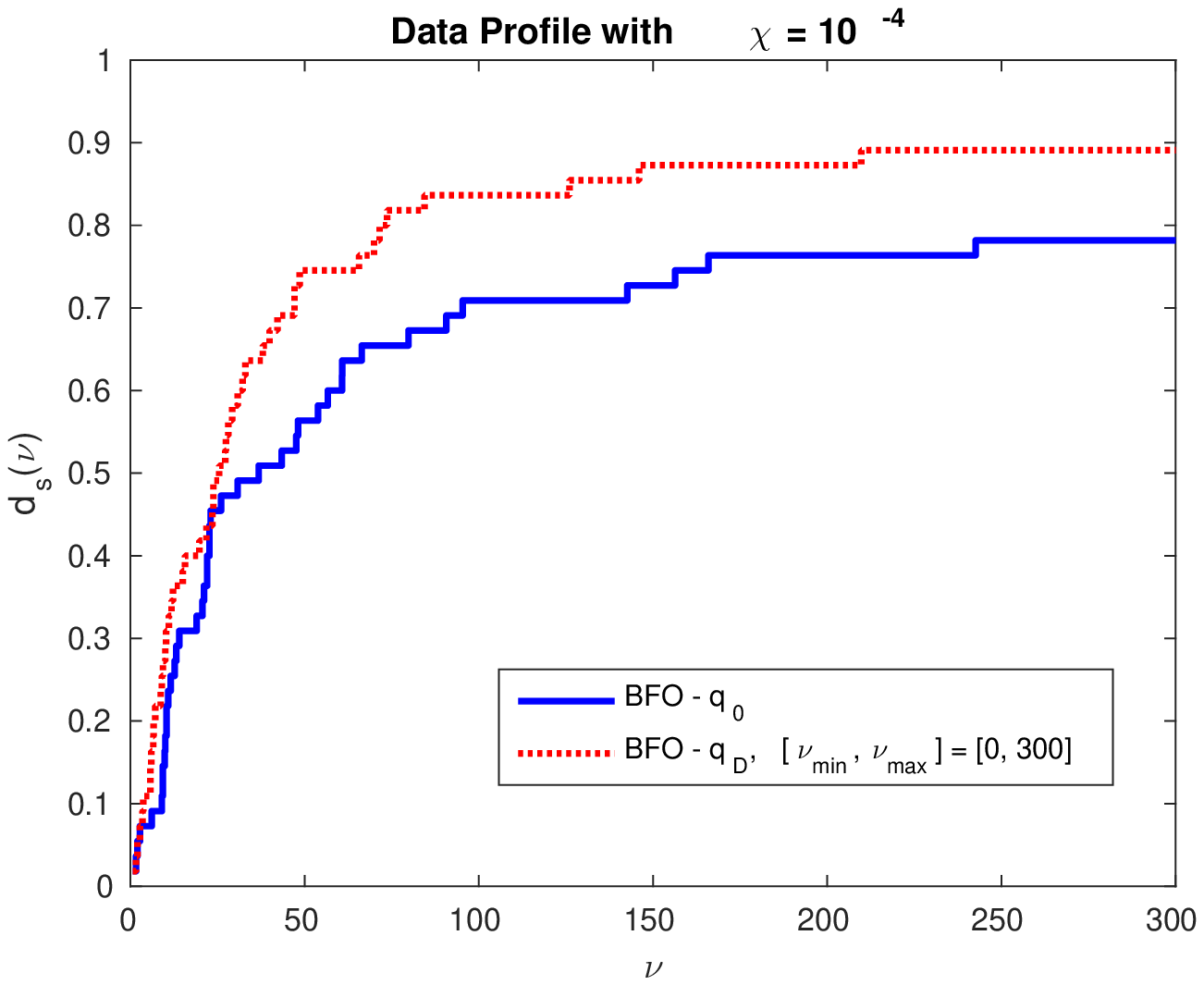}
  \caption{Zoom in the interval $[0, 300]$ of the data profiles obtained using BFO with
    $q_0$ and $q_{D}$ trained in the default  interval $[0, 2000]$ (left) and in the
    reduced interval $[0,300]$ (right).} 
  \label{figure:fig4}
\end{figure}

\begin{figure}[htb]
\centering
  \includegraphics[width=0.47\textwidth]{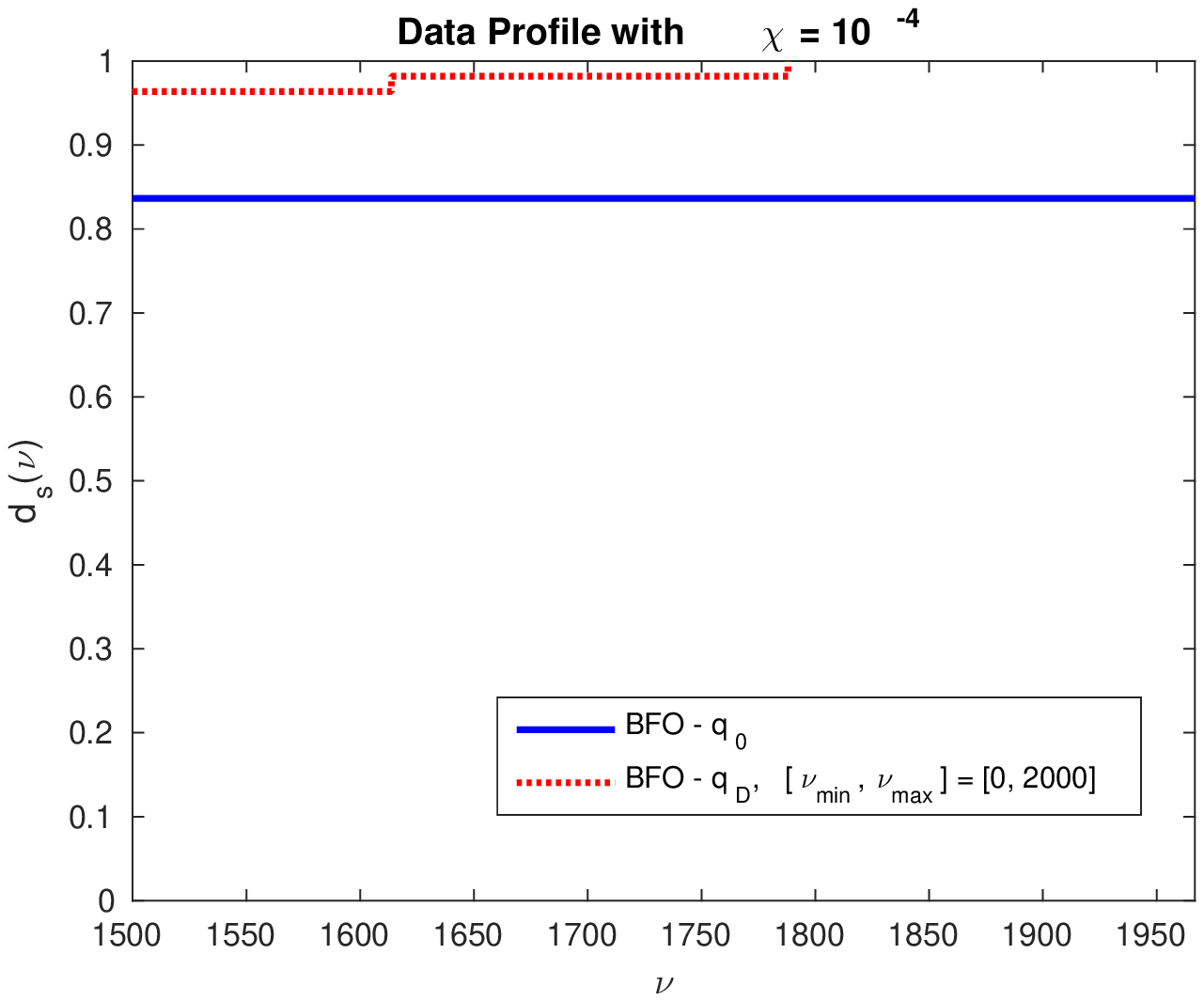}
  \includegraphics[width=0.47\textwidth]{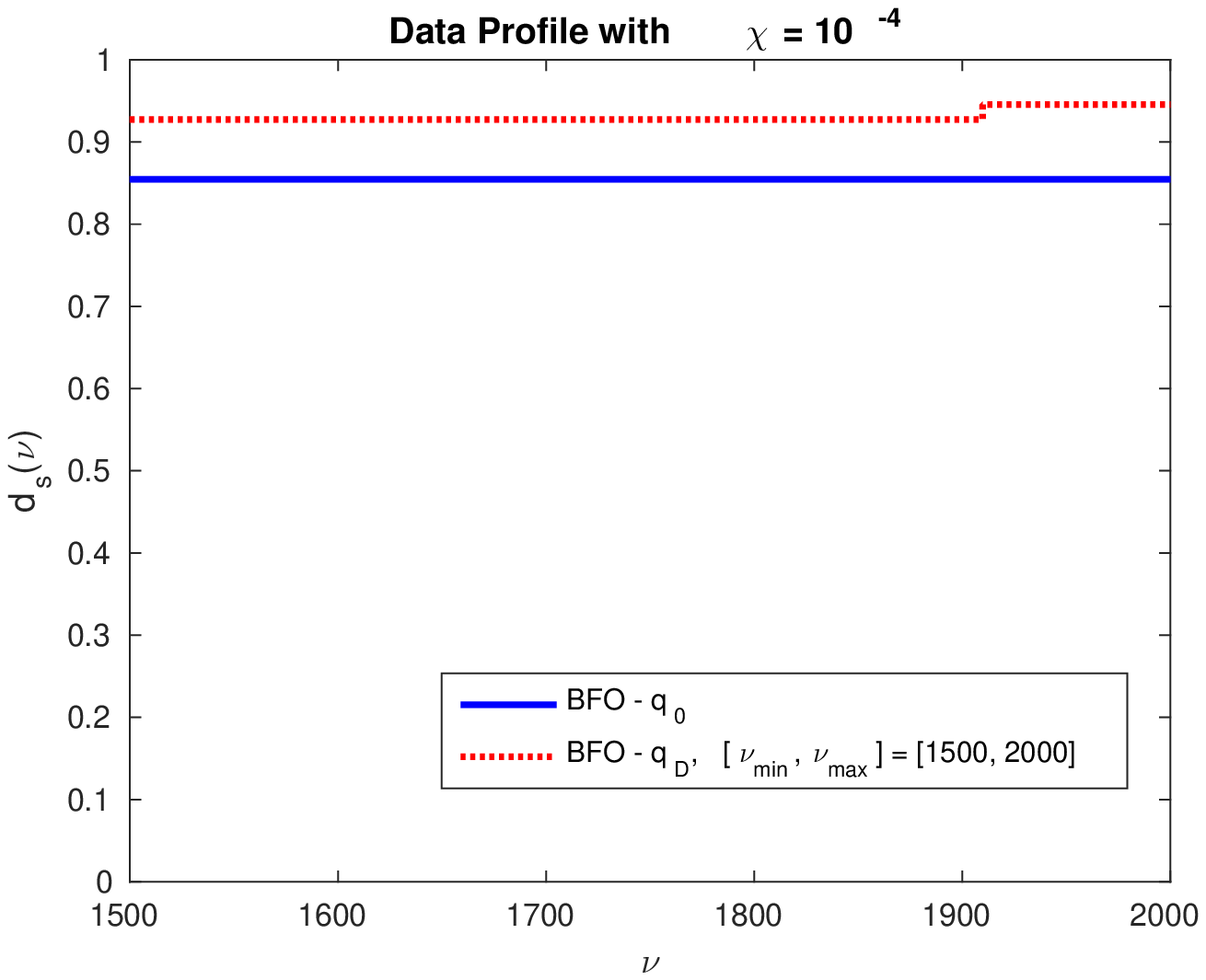}
  \caption{Zoom in the interval $[1500, 2000]$ of the data profiles obtained
    using BFO with $q_0$ and $q_{D}$ trained in the default  interval $[0,
      2000]$ (left) and in the reduced interval $[1500,2000]$ (right).} 
  \label{figure:fig5}
\end{figure}

We therefore repeated the training using the data-profile objective function
\req{DP} from the same initial parameter configuration $q_0$, but using
windows $[0, 300]$ and $[1500, 2000]$ instead of $[0,2000]$.  The resulting
%parameters are reported in Table~\ref{partuneDP} and the associated
profiles are presented in Figures~\ref{figure:fig4} and \ref{figure:fig5}. While the
expected improvement in efficiency using $[0,300]$ is clearly visible in the
first of these figures, the second shows that the procedure fails to produce
an improved reliability when using the window $[1500,2000]$, illustrating that
approximately and locally minimizing the training function ($\phi^D_\calP$ in
this case) does indeed sometimes produce sub-optimal solutions.

\section{Conclusion}

We have suggested how performance profiles and data profiles can be used to train
algorithms and have illustrated our proposal by an application to the BFO
package for derivative-free optimization. As expected, the results obtained
show that significant gains in performance are possible but not guaranteed.
The potential for improvement however suggests that the (careful) use of
the proposed techniques is a useful tool in algorithmic design.

\end{document}